\documentclass[11pt]{amsart}
\usepackage[latin1]{inputenc}
\usepackage{amsthm}
\usepackage{amssymb}
\usepackage{amsmath}
\usepackage{cancel}
\usepackage{a4wide}
\usepackage{verbatim}
\usepackage[normalem]{ulem}
\usepackage[
 colorlinks=true,urlcolor=blue,linkcolor=blue,citecolor=blue
]{hyperref}
\usepackage{hyperref}
\usepackage[capitalize]{cleveref}
\numberwithin{equation}{section} 

\theoremstyle{plain}
\newtheorem{thm}{Theorem}[section]

\theoremstyle{definition}
\newtheorem{defn}[thm]{Definition}
\theoremstyle{plain}
\newtheorem{prop}[thm]{Proposition}

\newtheorem{lem}[thm]{Lemma}
\newtheorem{cor}[thm]{Corollary}
\theoremstyle{remark}
\newtheorem{rem}[thm]{Remark}
\newtheorem{ex}[thm]{Example}

\usepackage{graphicx}
\usepackage{amscd}
\usepackage[all,knot,poly]{xy}
\usepackage[english]{babel}
\usepackage{xcolor}

\usepackage{mathrsfs}
\usepackage{hyperref}
\newcommand{\F}{\mathbb{F}}

\newcommand{\Z}{\mathbb{Z}}

\newcommand{\ones}{J}
\newcommand{\sJ}{\mathcal{J}}
\newcommand{\Rad}{{\rm Rad}}
\newcommand{\s}{\text{ss}}
\newcommand{\circulant}{\text{circ}}
\newcommand{\ord}{\text{ord}}

\newcommand{\diag}{\text{diag}}

\newcommand{\chebolu}[1]{{\color{teal}#1}}
\newcommand{\jon}[1]{{\color{purple}#1}}

\title[On the arithmetic of  join rings over finite fields]{On the arithmetic of  join rings \\ over finite fields}
\author[S.K. Chebolu, J. Merzel,J. Min\'a\v{c}, T. Nguyen,  N. D. Tan]{Sunil K. Chebolu, Jonathan L. Merzel, J\'an Min\'a\v{c}, \\ Tung T. Nguyen, Federico W. Pasini, Nguy$\tilde{\text{\^{E}}}$n Duy T\^{a}n}

\address{
Illinois State University}
\email{schebol@ilstu.edu}

\address{Soka University of America}
\email{jmerzel@soka.edu }

\address{
University of Western Ontario}
\email{minac@uwo.ca}

\address{Elmhurst University}
\email{tung.nguyen@elmhurst.edu}

\address{Huron University College}
\email{fpasini@uwo.ca}

\address{Faculty of Mathematics and Informatics, Hanoi University of Science and Technology}
\email{tan.nguyenduy@hust.edu.vn}

\date{}

\begin{document}
\thanks{Sunil Chebolu is partially supported by the Simons Foundation's Collaboration Grant for Mathematicians (516354). J\'an Min\'a\v{c} is partially supported by the Natural Sciences and Engineering Research Council of Canada (NSERC) grant R0370A01. J\'an Min\'a\v{c} also gratefully acknowledges 
Faculty of Sciences Distinguished Research Professorship award for 2020/21. J\'an Min\'a\v{c}, Tung T Nguyen, and Federico Pasini acknowledge the support of the Western Academy for Advanced Research. 
Nguy$\tilde{\text{\^{e}}}$n Duy T\^{a}n is funded by Vingroup Joint Stock Company and supported by Vingroup Innovation Foundation (VinIF) under the project code VINIF.2021.DA00030 and is partially supported by the Vietnam National
Foundation for Science and Technology Development (NAFOSTED) under grant number 101.04-2023.21}

\maketitle

\begin{abstract}

In this paper we consider some interesting and surprising interactions of several topics including representation theory, matrix algebra, and number theory. 
Given a collection $\{ G_i\}_{i=1}^d$ of finite groups and a ring $R$, we have previously introduced and studied certain foundational properties of the join ring $\sJ_{G_1, G_2, \ldots, G_d}(R)$.  This ring bridges two extreme worlds: matrix rings $M_n(R)$ on one end and group rings $R[G]$ on the other.  The construction of this ring was motivated by various problems in graph theory, network theory, nonlinear dynamics, and neuroscience.   In this largely self-contained paper, we continue our investigations of this ring, focusing more on its arithmetic properties.  We begin by constructing a generalized augmentation map that gives a structural decomposition of this ring. This decomposition allows us to compute the zeta function of the join of group rings. We show that the join of group rings is a natural home for studying the concept of simultaneous primitive roots for a given set of primes. This concept is related to the order of the unit group of the join of group rings. Finally, we characterize the join of group rings over finite fields with the property that the order of every unit divides a fixed number. Remarkably, Mersenne and Fermat primes unexpectedly emerge within the context of this exploration.
\end{abstract}
\tableofcontents
\section{Introduction}
Let $G$ be a finite group. The concept of $G$-circulant matrices, defined in \ref{defn:circulant}, has a rich mathematical history. Dedekind initially introduced these matrices during his study of normal bases for Galois extensions. His focus was on understanding the factorizations of the determinants of $G$-circulant matrices. While his success was notable in cases where $G$ is abelian, his progress in the non-abelian realm was limited. This work led to correspondence with Frobenius in 1896. Subsequently, Frobenius made a pivotal discovery. He showed that the determinant of a generic $G$-circulant matrix decomposes into a product of irreducible factors over the field of complex numbers corresponding to the linear irreducible representations of the group $G$.
 In particular, when $G$ is a cyclic group, we have an explicit description of the spectrum of $G$-circulant matrices. This description is often referred to as the Circulant Diagonalization Theorem in the literature (see \cite{davis2013circulant} for an extensive treatment of this topic). Due to their elegance and explicit nature, circulant matrices have found applications in many scientific fields, such as spectral graph theory, coding theory, neuroscience, and nonlinear dynamics (see \cite{cir1, ko3, davis2013circulant, CM1,cir2, kanemitsu2013matrices, cir3, nguyen2023equilibria, townsend2020dense}). In particular, in \cite{ko3}, using the spectral decomposition of a circulant network, we are able to explain various traveling wave patterns in networks of phase oscillators.

 In \cite{CM2}, we introduce a natural generalization of $G$-circulant matrices. More precisely, given a collection of finite groups $G_1, G_2, \ldots, G_d$ and a ring $R$, we introduce the join ring $\sJ_{G_1, G_2, \ldots, G_d}(R)$ (see Section \ref{subsection:join_ring} for the precise definition of this ring). When $d=1$, the ring $\sJ_{G}(R)$ is exactly the ring of all $G$-circulant matrices with entries in $R$. Furthermore, $\sJ_{G}(R)$ is naturally isomorphic to the group ring $R[G].$ We also remark that when all $G_i$ are the trivial group, the join ring is naturally isomorphic to $M_d(R)$, the ring of all square matrices of size $d \times d$ with coefficients in $R$. The introduction of the join ring $\sJ_{G_1, G_2, \ldots, G_d}(R)$ is motivated by a construction in graph theory known as the joined union of graphs, and by a desire to understand nonlinear dynamics in multilayer networks of oscillators (see \cite{CM1, CM1_b, nguyen2023equilibria,nguyen2023broadcasting}). In \cite{CM2}, we discuss some fundamental ring-theoretic properties of the join ring $\sJ_{G_1, G_2, \ldots, G_d}(R)$ such as its center, its semisimplicity, its Jacobson radical, the structure of its unit group, and much more. In this article, we discuss some further properties of this ring, focusing on the case that $R$ is a finite field. This article presents our continuing effort to develop a systematic understanding of the join ring $\sJ_{G_1, G_2, \ldots,G_d}(R).$ 
 We have made a concerted effort to ensure that our work is accessible to a broad readership. To this end, we provide a self-contained review of the fundamental notions and key results required for a complete understanding of the text.

We now summarize our main results. The definitions of the join of group rings and the associated zeta functions can be found in \cref{sec:ring_theoretic} and Section~\ref{sec:zeta}, respectively. Our first result is a structural decomposition of the join rings.

\begin{thm}[Decomposition of Join Rings]
Let $G_1, \ldots, G_d$ be finite groups and $H_i \trianglelefteq G_i$ such that $|H_i|$ is invertible in a unital ring $R$. Then, there exists a ring isomorphism
\[
J_{G_1, \ldots, G_d}(R) \cong J_{G_1/H_1, \ldots, G_d/H_d}(R) \times \prod_{i=1}^d \Delta_R(G_i, H_i),
\]
where $\Delta_R(G_i, H_i)$ is the kernel of the augmentation map $R[G_i] \to R[G_i/H_i]$.
\end{thm}

Note that, in the special case when all the   $|G_i|$'s are invertible in $R$, we get
\[ \sJ_{G_1, G_2, \ldots, G_d}(R) \cong M_{d}(R) \times \prod_{i=1}^d \Delta_R(G_i).\] 

The above structural decomposition helps compute the zeta function for the join of group rings. We refer the reader to \cref{sec:ring_theoretic} for the definition of these rings and to \cref{sec:zeta} for the definition of their zeta functions.

\begin{thm}[Zeta Function of Join Rings]

Let $\F_q$ be a finite field, and suppose that $|G_i|$ is invertible in $\F_q$ for all $1 \le i \le d$. Then the zeta function of the join ring satisfies
\[
\zeta_{J_{G_1, \ldots, G_d}(\F_q)}(s) = (1 - q^{-s})^{d - 1} \prod_{i=1}^d \zeta_{\F_q[G_i]}(s).
\]
\end{thm}

For the general case when some $|G_i|$ is not invertible in $\F_q$, we refer the reader to  \cref{cor:general-zeta}. We use these results to explicitly compute the zeta function in a number of examples; see \cref{cor:zeta_abelian_groups} and various other examples in Section \ref{sec:zeta}.

We then use these zeta functions to study $q$-rooted primes. A prime $p$ is said to be $q$-rooted if $q$ is a primitive root modulo $p$. $2$-rooted primes and their characterizations were studied in  \cite{chebolu2015characterizations}. In the current paper, we extend those results to the odd primary case and also study simultaneous  $q$-rooted primes in the framework of zeta functions.

\begin{thm}[Characterization of $q$-Rooted Primes]
Let $q, p_1, \ldots, p_d$ be prime numbers with $p_i \ne q$ for all $i$. The following are equivalent:
\begin{enumerate}
    \item Each $p_i$ is a $q$-rooted prime (i.e., $q$ is a primitive root modulo $p_i$).
    \item The order of the pole at $s = 0$ of $\zeta_{J_{\mathbb{Z}/p_1, \ldots, \mathbb{Z}/p_d}(\F_q)}(s)$ is $d+1$.
    \item The order of the unit group of $J_{\mathbb{Z}/p_1, \ldots, \mathbb{Z}/p_d}(\F_q)$ is
    \[
    \prod_{i=1}^d (q^{p_i - 1} - 1) \cdot \prod_{i = 0}^{d - 1}(q^d - q^i).
    \]
\end{enumerate}
\end{thm}

We end the paper with an investigation of  $\Delta_{p^r}$-rings. A ring is said to be a  $\Delta_{p^r}$-ring if $u^{p^r} = 1$ for all units $u$ in the ring. These rings were introduced 
in  \cite{chebolu2015characterizations}, where the following question was raised: When is the group algebra $kG$ a $\Delta_p$-ring? In  \cite{chebolu2015characterizations}, the authors addressed this for the case when $G$ is an abelian group and $r=1$.  Here we extend those results to all finite groups and $r \ge 1$ (\ref{thm:delta_classification}), and also to  join algebras defined over a  finite field. 

\begin{thm}[ Classification of Join Rings that are $\Delta_{p^r}$-Rings]
Let $d \ge 2$. Then the join ring $J_{G_1, \ldots, G_d}(\F_q)$ is a $\Delta_{p^r}$-ring if and only if:
\begin{enumerate}
    \item $p = q = 2$,
    \item Each $G_i$ is a $2$-group,
    \item At most one $G_i$ is trivial.
    \item $2^r \ge \max_{1 \le i \le d} \exp(U_1(\F_2[G_i]))$, where $U_1(\F_2[G_i])$ is the set of normalized units in $\F_2[G_i]$.
\end{enumerate}
\end{thm}

\subsection{Outline}
The structure of this article is as follows. In \cref{sec:ring_theoretic}, we study some further ring-theoretic properties of the join ring $\sJ_{G_1, G_2, \ldots, G_d}(R)$. Among various things that we discover, we discuss a natural construction of the generalized augmentation map (a special case of this construction is previously discussed in \cite{CM2}). \cref{sec:zeta} studies the zeta functions of the join ring $\sJ_{G_1, G_2, \ldots, G_d}(R)$ when $R=\mathbb{F}_q$ is a finite field. More precisely, we describe how to explicitly calculate the zeta function of $\sJ_{G_1, G_2, \ldots, G_d}(\F_q)$ in terms of the zeta functions of 
$\sJ_{G_i}(\F_q)$. In \cref{sec:rooted_prime}, we discuss the order of the unit group of $\sJ_{G_1, G_2, \ldots, G_d}(\F_q)$ and explain its connection with Artin's conjecture on primitive roots. For instance, we find several equivalent conditions for when a given prime $q$ is simultaneously a primitive root for a set of primes $\{p_i \}$. These equivalent conditions are based on the cardinality of the unit group and the order of the pole of the zeta function for the join ring $\sJ_{\Z/p_1, \Z/p_2, \ldots, \Z/p_d}(\F_q)$.
Finally, in \cref{sec:delta_ring}, we classify all join rings $\sJ_{G_1, G_2, \ldots, G_d}(\F_q)$ that have the property that every unit $u$ in them satisfies $u^{p^r}=1$, where $p$ is a prime number and $r$ is a positive integer.  Such rings are called $\Delta_{p^r}$ rings and they are well-studied in the literature; see \cite{chebolu2012special,  chebolu2015characterizations, chebolu2013special}.
In particular, \cite{chebolu2015characterizations} focuses on the aforementioned property when $d=1$, $r=1$, and the finite group involved is abelian. The results of Section \ref{sec:delta_ring}  advance beyond these parameters, broadening the outcome to encompass all finite groups and positive integer values of $d$ and $r$.

\section*{Acknowledgements}
We thank Professors Kazuya Kato, Tsit Yuen Lam, and Michel Waldschmidt for their helpful correspondence and encouragement. We also thank the referee for the comments and suggestions that helped us to polish some of our exposition.

\section{Some ring-theoretic properties of $\sJ_{G_1, G_2, \ldots, G_d}(R)$.} 
\label{sec:ring_theoretic}
\subsection{The ring of $G$-circulant matrices}
Let $G = \{g_1 (=e), g_2, \ldots, g_n\}$ be a finite group of order $n$ (note that we have fixed an ordering on $G$). We first recall the definition of a $G$-circulant matrix (for more details, see \cite{CM2, hurley2006, kanemitsu2013matrices}). 
\begin{defn} \label{defn:circulant}
An $n \times n$ $G$-circulant matrix over $R$ is an $n\times n$ matrix 
\[
A=\left[ 
\begin{array}{llll}
a_{g_{1},g_{1}} & a_{g_{1},g_{2}} & \cdots  & a_{g_{1},g_{n}} \\ 
a_{g_{2},g_{1}} & a_{g_{2},g_{2}} & \cdots  & a_{g_{2},g_{n}} \\ 
\vdots  & \vdots  & \ddots  & \vdots  \\ 
a_{g_{n},g_{1}} & a_{g_{n},g_{2}} & \cdots  & a_{g_{n},g_{n}}%
\end{array}%
\right] 
\]%
over $R$ with the property that
for all $g,g_i,g_j \in G,$ $a_{g_i,g_j}=a_{gg_i,gg_j}$.\
\end{defn}
We remark that a $G$-circulant matrix $A$ is completely determined by its first row and the multiplication table of $G$, as we must have $a_{g_i,g_j}=a_{g_1,g_i^{-1}g_j}$. For simplicity, we sometimes write $A=\circulant([a_g]_{g \in G})$ (where we define the doubly indexed quantity $a_{g_i,g_j}$ by $a_{g_i^{-1}g_j}$).  Let $\sJ_{G}(R)$ be the set of all $G$-circulant matrices over $R$ and   
\[ R[G]=\left\{\sum_{g \in G} a_g g \big|a_g \in R \right \} ,\] 
 the group ring of $G$ with coefficients in $R.$ In \cite{CM2}, we reproved the following theorem of Hurley. 
\begin{prop} (Hurley) \label{prop:hurley}
The map $\alpha: R[G] \to \sJ_G(R)$ sending 
\[ \sum_{g \in G} a_g g \mapsto \circulant([a_g]_{g \in G}), \]
is a ring isomorphism. In particular, under this isomorphism, units in the group ring correspond to invertible $G$-circulant matrices.
\end{prop}

\subsection{The join ring $\sJ_{G_1, G_2, \ldots, G_d}(R)$} \label{subsection:join_ring}
We recall the definition of the join matrix (see \cite[Definition 3.1]{CM2}).
\begin{defn}
Let $R$ be a (unital, associative) ring, $G_{1},\ldots ,G_{d}$  finite
groups of respective orders $k_{1},\ldots ,k_{d}$, and let $C_{i}$ be $G_{i}$%
-circulant ($1\leq i\leq d$) over $R$. \ By a join of $C_{1},\ldots ,C_{d}$
over $R$, we mean a matrix of the form 

\begin{equation}
\label{eq:join circulant matrix}\tag{$\ast$}
A=\begin{bmatrix} 
C_1&a_{12} J_{k_1,k_2} &\cdots & a_{1d}J_{k_1,k_d}\\
a_{21}J_{k_2,k_1} &C_2 &\cdots & a_{2d}J_{k_2,k_d}\\
\vdots&\vdots& &\vdots\\
a_{d1}J_{k_d,k_1}&a_{d2}J_{k_d,k_2}&\cdots& C_d
\end{bmatrix}
, 
\end{equation}
where $a_{ij}\in R\ (1\leq i\neq j\leq d)$ and $J_{r,s}$ denotes the $%
r\times s$ matrix, all of whose entries are $1\in R$.  
\end{defn} 
We remark that we came upon the concept of a join matrix through our work on multilayer networks of phase oscillators (see \cite{CM1, CM1_b, nguyen2023broadcasting}). As in \cite{CM2}, we will denote by $\sJ_{G_{1},\ldots,G_{d}}(R)$, the set of all such joins as the $C_{i\text{ }}$vary independently through all $G_{i}$-circulant matrices ($1\leq i\leq d$) and
the $a_{ij}$ vary independently through all elements of $R$ ($1\leq i\neq j\leq d$). In \cite{CM2}, we showed the following.
\begin{prop} (\cite[Section 3]{CM2}
$\mathcal{J}_{G_1,\dots, G_d}(R)$ has the structure of a unital ring. 
Furthermore, there is an augmentation map $\epsilon \colon \mathcal{J}_{G_1,\dots, G_d}(R) \rightarrow M_d(R)$ that generalizes the augmentation map on group rings.
\end{prop}

Since we identified $\mathcal{J}_{G_1,\dots, G_d}(R)$ as a subring of a matrix ring over $R$, it is clear that in the case when $R$ is a field $k$, $\mathcal{J}_{G_1,\dots, G_d}(k)$ has the structure of a $k$-algebra.

\subsection{The generalized augmentation map} 
Let $G$ be a finite group and $H$ be a normal subgroup of $G$. Then, there is a canonical ring map known as the augmentation map
\begin{equation} \label{eq:classical_aug}
\epsilon: R[G] \to R[G/H],
\end{equation}
which extends the quotient map $G \rightarrow G/H$ that sends $g \mapsto \bar{g}$. When $H=G$, this is exactly the standard augmentation map $\epsilon: R[G] \to R$ mentioned in the previous section. More concretely, this augmentation map is defined by  
\[ \epsilon \left(\sum_{g \in G} a_g g \right) = \sum_{g \in G} a_g .\] 
In this section, we show that there is a natural analog of this augmentation map in the setting of the join ring $\sJ_{G_1, G_2, \ldots, G_d}(R).$ More precisely, let $G_i$ be a finite group  and $H_i$ a normal subgroup of $G_i$ for all $1 \leq i \leq d$. Suppose that the orders of $G_i, H_i, G_i/H_i$ are $k_i, r_i, s_i$ respectively (so $k_i=r_i s_i$). Let us consider the following map
\begin{equation} \label{eq:generalized_aug}
\epsilon: \sJ_{G_1, G_2, \ldots, G_d}(R) \to \sJ_{G_1/H_1, G_2/H_2, \ldots, G_d/H_d}(R), 
\end{equation}
defined by 
\begin{equation*}
\label{eq:augmentation}
\begin{bmatrix} 
C_1&a_{12} J_{k_1,k_2} &\cdots & a_{1d}J_{k_1,k_d}\\
a_{21}J_{k_2,k_1} &C_2 &\cdots & a_{2d}J_{k_2,k_d}\\
\vdots&\vdots& &\vdots\\
a_{d1}J_{k_d,k_1}&a_{d2}J_{k_d,k_2}&\cdots& C_d
\end{bmatrix}  \mapsto \begin{bmatrix} 
\epsilon(C_1)& r_2 a_{12} J_{s_1,s_2} &\cdots & r_d a_{1d}J_{s_1,s_d}\\
r_1 a_{21}J_{s_2,s_1} &\epsilon(C_2) &\cdots & r_d a_{2d}J_{s_2,s_d}\\
\vdots&\vdots& &\vdots\\
r_1 a_{d1}J_{s_d,s_1}& r_2 a_{d2}\sJ_{s_d,s_2}&\cdots& \epsilon(C_d).
\end{bmatrix} 
\end{equation*} 
Here $\epsilon$ is the classical augmentation map $R[G_i] \to R[G_i/H_i]$ as defined in Equation \ref{eq:classical_aug}. We remark that the row and column sum of a $G$-circulant matrix $A=\text{circ}([a_g]_{g \in G})$ are all equal to $\sum_{g \in G} a_g$. This type of matrix has a special name which we now recall. 
\begin{defn} \label{defn:semimagic}(see \cite{murase1957semimagic})
Let $R$ be a ring. A matrix $A \in M_n(R)$ is called a semimagic square if its row and column sums are equal; i.e., there exists a constant $\sigma(A)$ such that 
\[ \sum_{i=1}^n a_{ij} = \sum_{j=1}^n a_{ij} = \sigma(A). \]

\end{defn}

We have the following.

\begin{prop} \label{prop:ring_hom_26}
The map $\epsilon: \sJ_{G_1, G_2, \ldots, G_d}(R) \to \sJ_{G_1/H_1, G_2/H_2, \ldots, G_d/H_d}(R)$ is a ring homomorphism. 
\end{prop}
\begin{proof}
This follows from direct calculations. Two key identities are the following.
\begin{enumerate}
    \item $J_{m,n} \times J_{n,p}=n J_{m,p}.$
    \item $A J_{m,n}= \epsilon(A) J_{m,n}$ where $A$ is a semimagic square of size $m \times m$ and $\epsilon(A)$ is the row sum of $A.$ Similarly $J_{m,n}B = \epsilon(B) J_{m,n}$ if $B$ is a semimagic square of size $n \times n$.
\end{enumerate}
\end{proof} 

Given a group homomorphism $\varphi :G\rightarrow G^{\prime }$ we have
functorially a map $\Phi ^{G,G^{\prime }}:R[G]\rightarrow R[G^{\prime }]$. \
(We could, of course, also permit a homomorphism from $R$ to another ring $%
R^{\prime }$.) In the case where $H\vartriangleleft G$ and $\varphi $ is the
canonical map $\pi _{G,G/H}:G\rightarrow G/H,$ the map $\Phi ^{G,G/H}$ is
what we above called the augmentation map $\epsilon $, which we will
temporarily designate $\epsilon ^{G,G/H}.$ \ By functoriality, we mean
that if we also have $\varphi ^{\prime }:G^{\prime }\rightarrow G^{\prime
\prime }$ then $\Phi ^{G,G^{\prime \prime }}=\Phi ^{G^{\prime },G^{\prime
\prime }}\circ \Phi ^{G,G^{\prime }}$. \ We therefore have

\begin{lem}
(a) \ If $\varphi :G\rightarrow G^{\prime }$ is a homomorphism and if $%
H\vartriangleleft G$ with $\varphi (H)<H^{\prime }\vartriangleleft G^{\prime
}$, then we have a commutative diagram%
\[
\begin{array}{ccc}
R[G] & \overset{\epsilon ^{G,G/H}}{\rightarrow } & R[G/H] \\ 
\Phi ^{G,G^{\prime }}\downarrow \ \ \ \ \ \ \  &  & \ \ \ \ \ \ \ \ \ \ \
\downarrow \Phi ^{G/H,G^{\prime }/H^{\prime }} \\ 
R[G^{\prime }] & \overset{\epsilon ^{G^{\prime },G^{\prime }/H^{\prime }}}%
{\rightarrow } & R[G^{\prime }/H^{\prime }]%
\end{array}%
\]%
(b) \ If $H,K\vartriangleleft G$ with $H<K$ then $\epsilon
^{G,G/K}=\epsilon ^{G/H,G/K}\circ \epsilon ^{G,G/H}.$
\end{lem}

\begin{proof}
(a) \ As $\epsilon ^{G,G/H}=\Phi ^{G,G/H}$ and $\epsilon ^{G^{\prime
},G^{\prime }/H^{\prime }}=\Phi ^{G^{\prime },G^{\prime }/H^{\prime }}$, by
functoriality both compositions in the diagram  are $\Phi ^{G,G^{\prime
}/H^{\prime }}$ induced by $\pi _{G^{\prime },G^{\prime }/H^{\prime }}\circ
\varphi $.\newline
(b) Apply (a) to the case $G^{\prime }=G,\ H^{\prime }=K.$
\end{proof}

We would like to generalize the above to join rings. \ Unfortunately, given
homomorphisms $\varphi _{i}:G_{i}\rightarrow G_{i}^{\prime },\ i=1,\cdots ,d$
there is in general no apparent naturally induced map $\Phi :J_{G_{1},\cdots
,G_{d}}(R)\rightarrow J_{G_{1}^{\prime },\cdots ,G_{d}^{\prime }}(R)$. \ The
problem is that the natural image of $J_{k_{i},k_{i}}\in R[G_{i}]$ is not a
multiple of $J_{k_{i}^{\prime },k_{i}^{\prime }}\in R[G_i^{\prime }]$;
however, this holds when the maps $\varphi _{i}:G_{i}\rightarrow
G_{i}^{\prime }$ are surjective. \ In that case, we are up to isomorphism
back to the earlier situation of \cref{prop:ring_hom_26}, and we will also denote the
map $\epsilon $ of that proposition by $\epsilon
^{\{G_{i}\},\{G_{i}^{\prime }\}}$. \ Corresponding to the lemma above, we
now have

\begin{prop}
(a) \ Given surjective homomorphisms $\varphi _{i}:G_{i}\rightarrow
G_{i}^{\prime },\ i=1,\cdots ,d$  and normal subgroups $H_{i}%
\vartriangleleft G_{i}$ with $\varphi (H_{i})<H_{i}^{\prime
}\vartriangleleft G_{i}^{\prime }$, we then have a commutative diagram 
\[
\begin{array}{ccc}
J_{G_{1},\cdots ,G_{d}}(R) & \overset{\epsilon
^{\{G_{i}\},\{G_{i}/H_{i}\}}}{\rightarrow } & J_{G_{1}/H_{1},\cdots
,G_{d}/H_{d}}(R) \\ 
\epsilon ^{\{G_{i}\},\{G_{i}^{\prime }\}}\downarrow \ \ \ \ \ \ \  &  & \
\ \ \ \ \ \ \ \ \ \ \downarrow \epsilon ^{\{G_{i}/H_{i}\},\{G_{i}^{\prime
}/H_{i}^{\prime }\}} \\ 
J_{G_{1}^{\prime },\cdots ,G_{d}^{\prime }}(R) & \overset{\epsilon
^{\{G_{i}^{\prime }\},\{G_{i}^{\prime }/H_{i}^{\prime }\}}}{\rightarrow } & 
J_{G_{1}^{\prime }/H_{1}^{\prime },\cdots ,G_{d}^{\prime }/H_{d}^{\prime
}}(R)%
\end{array}%
\]%
(b) \ If $H_{i},K_{i}\vartriangleleft G_{i}$ with $H_{i}<K_{i},i=1,\cdots
,d\ $ then 
\[
\epsilon ^{\{G_{i}\},\{G_{i}/K_{i}\}}=\epsilon
^{\{G_{i}/H_{i}\},\{G_{i}/K_{i}\}}\circ \epsilon
^{\{G_{i}\},\{G_{i}/H_{i}\}}.
\]
\end{prop}

\begin{proof}
(a) \ For a matrix%
\[
\left( 
\begin{array}{llll}
C_{1} & a_{12}J_{k_{1},k_{2}} & \cdots  & a_{1d}J_{k_{1},k_{d}} \\ 
a_{21}J_{k_{2},k_{1}} & C_{2} & \cdots  & \alpha _{2d}J_{k_{2},k_{d}} \\ 
\vdots  & \vdots  & \ddots  & \vdots  \\ 
\alpha _{d1}J_{k_{d},k_{1}} & \alpha _{d2}J_{k_{d},k_{2}} & \cdots  & C_{d}%
\end{array}%
\right) \in J_{G_{1},\cdots ,G_{d}}(R)
\]%
we consider the image of its blocks under the compositions $\epsilon
^{\{G_{i}/H_{i}\},\{G_{i}^{\prime }/H_{i}^{\prime }\}}\circ \epsilon
^{\{G_{i}\},\{G_{i}/H_{i}\}}$ and $\epsilon ^{\{G_{i}^{\prime
}\},\{G_{i}^{\prime }/H_{i}^{\prime }\}}\circ \epsilon
^{\{G_{i}\},\{G_{i}^{\prime }\}}$. \ For the blocks on the diagonal, these
images coincide by part (a) of the above lemma. \ The $i,j$ block with $%
i\neq j$ accumulates factors of $\left\vert H_{j}\right\vert $ and $%
\left\vert (G_{j}/H_{j})/(G_{j}^{\prime }/H_{j}^{\prime })\right\vert $
under the composition $\epsilon ^{\{G_{i}/H_{i}\},\{G_{i}^{\prime
}/H_{i}^{\prime }\}}\circ \epsilon ^{\{G_{i}\},\{G_{i}/H_{i}\}}$ and
factors of $\left\vert G_{j}\right\vert /\left\vert G_{j}^{\prime
}\right\vert $ and $\left\vert H_{j}^{\prime }\right\vert $ under the
composition $\epsilon ^{\{G_{i}^{\prime }\},\{G_{i}^{\prime
}/H_{i}^{\prime }\}}\circ \epsilon ^{\{G_{i}\},\{G_{i}^{\prime }\}}$, so
the image either way is $\frac{\left\vert G_{j}\right\vert \left\vert
H_{j}^{\prime }\right\vert }{\left\vert G_{j}^{\prime }\right\vert }%
a_{ij}J_{G_{i}^{\prime }/H_{i}^{\prime },G_{j}^{\prime }/H_{j}^{\prime }}$.%
\newline
(b) \ Apply part (a) to the case $G_{i}^{\prime }=G_{i},\ \ H_{i}^{\prime
}=K_{i}$.
\end{proof}

\subsection{A decomposition of $\sJ_{G_1, G_2, \ldots, G_d}(R)$.}
Let $\Delta_R(G,H) := ker(R[G] \to R[G/H])$ be the kernel of the augmentation map as defined in Equation \ref{eq:classical_aug} (when $R$ is clear from the context, we will simply write $\Delta(G,H)$.) Suppose further that $|H|$ is invertible in $R$. Let 
\[ e_{H}=\frac{1}{|H|} \sum_{h \in H} h .\] 
It can be shown that $e_H$ is a central idempotent in $R[G]$; see also \cite[Lemma 3.6.6]{milies2002introduction}. Furthermore, by \cite[Proposition 3.6.7]{milies2002introduction}, we have 
\begin{prop} \label{prop:splitting}
We have a direct product of rings 
\[ R[G] \cong R[G]e_H \times R[G](1-e_H) .\] 
Furthermore 
\[ R[G] e_H \cong R[G/H], \]
and 
\[ R[G](1-e_H) = \Delta_R(G,H) .\] 
\end{prop}

\begin{cor}(see \cite[Corollary 3.6.9]{milies2002introduction})
Suppose that $|G|$ is invertible in $R$. Let $\Delta_R(G)$ be the augmentation ideal. Then 
\[ R[G] \cong R \times \Delta_R(G). \] 
\end{cor}

We can generalize this proposition to the join ring as follows. 
\begin{thm} \label{thm:main}
Let $G_1, \ldots, G_d$ be finite groups. For $1 \leq i \leq d$, let $H_i$ be a normal subgroup such that $|H_i|$ is invertible in $R$. Then, there exists an isomorphism 
\[ \sJ_{G_1, \ldots, G_d}(R) \cong \sJ_{G_1/H_1, \ldots, G_d/H_d}(R) \times \prod_{i=1}^d \Delta_R(G_i,H_i) .\] 

\end{thm}
\begin{proof}
Let $f_i=f_{G_i}=1-e_{H_i} \in R[G_i]$ where $e_{H_i}$ is defined as above. Since the ring of all $G_i$-circulant matrices is isomorphic to the group ring $R[G_i]$, we can also consider $f_i$ as a $G_i$-circulant matrix. Let $\tilde{f}_i$ be the following matrix in $\sJ_{G_1, G_2, \ldots, G_d}(R)$
\[
\tilde{f}_i =\left[\begin{array}{c|c|c|c}
0 & 0 & \dots & 0 \\
\hline
0 & 0 & \dots & 0 \\
\hline
\vdots & \vdots & f_i & \vdots\\
\hline
0 & 0 & \dots & 0
\end{array}
\right].
\]
In other words, all blocks of $\tilde{f}_i$, except the $i$-diagonal block which is $f_i$,  are $0.$  Additionally, we define 
\[ \tilde{f}_{d+1}=I_n-\sum_{i=1}^d \tilde{f}_i=\bigoplus_{i=1}^{d} e_{G_i}. \] 
Then we have the following ring isomorphism 
\[ \sJ_{G_1, \ldots, G_d}(R) \cong \tilde{f}_{d+1} \sJ_{G_1, \ldots, G_d}(R) \times \prod_{i=1}^d \tilde{f}_{i} \sJ_{G_1, \ldots, G_d}(R) .\] 
We can see that for $1 \leq i \leq d$ 
\[ \tilde{f}_{i} \sJ_{G_1, \ldots, G_d}(R) \cong \Delta_R(G_i, H_i) .\] 

Additionally, the augmentation map 
\[ \epsilon: \sJ_{G_1, G_2, \ldots, G_d}(R) \to \sJ_{G_1/H_1, \ldots, G_d/H_d}(R) \] 
induces a ring isomorphism 
\[ \epsilon: \tilde{f}_{d+1} \sJ_{G_1, G_2, \ldots, G_d}(R) \to \sJ_{G_1/H_1, \ldots, G_d/H_d}(R).\] 
\end{proof}

Here is a direct corollary of this theorem.
\begin{cor}(See also \cite[Theorem 3.16]{CM2}) \label{cor:decomposition}
Suppose that $|G_i|$ are invertible in $R$. Then 
\[ \sJ_{G_1, G_2, \ldots, G_d}(R) \cong M_{d}(R) \times \prod_{i=1}^d \Delta_R(G_i) .\] 

\end{cor}
\section{Zeta function of the join ring $\sJ_{G_1, G_2, \ldots, G_d}(\F_q)$} 
\label{sec:zeta}
Let $\F_q$ be the finite field with $q=p^r$ elements where $p$ is a prime number. In this section, we study the zeta function of the join ring $\sJ_{G_1, G_2, \ldots, G_d}(\F_q).$ We first recall the definition of the zeta function of a finite-dimensional $\F_q$-algebra as defined in \cite{fukaya1998hasse}.

First, consider the case where $R$ is a commutative finite dimensional $\F_q$-algebra. The Hasse-Weil zeta function of $R$ is defined to be 
\begin{equation} \label{eq:hasse_weil}
\zeta_{R}(s)=\prod_{m \subset R} (1-\#(R/m)^{-s})^{-1} .
\end{equation}
where $m$ runs over all maximal ideals of $R$ (see \cite{fukaya1998hasse}). As observed in \cite{fukaya1998hasse}, when $R$ is not commutative, the Hasse-Weil zeta function of $R$ can be defined as follows. (We refer readers to \cite{kurokawa1984some, kurokawa1989special} for some further motivations for this definition.)
\begin{defn}(see \cite{fukaya1998hasse})
Let $R$ be a 
finite-dimensional $\F_q$-algebra. The following Euler product gives the Hasse-Weil zeta function of $R$ 
\begin{equation} \zeta_{R}(s) = \prod_{M}(1-  |\text{End}_{R}(M)|^{-s})^{-1},
\end{equation}
where $M$ runs over the isomorphism classes of (finite) simple left $R$-modules. 
\end{defn}
We remark that since $R$ is a finite ring, all simple left $R$-modules are automatically finite. Furthermore, by \cite[Lemma 2.7.1]{fukaya1998hasse}, another equivalent definition of $\zeta_{R}(s)$ is 
\begin{equation} \zeta_{R}(s) = \prod_{m \in \mathfrak{P}(R)} (1-N(\frak{m})^{-s})^{-1} ,
\end{equation}
where $\mathfrak{P}(R)$ is the set of all two-sided ideals $\frak{m}$ in $R$ such that $R/\frak{m}$  is isomorphic to a matrix ring $M_r(k)$ with $k$  a finite extension of $\F_q$ and $N(\frak{m})= |k|.$

For a finite-dimensional $\F_q$-algebra $R$, we denote 
\begin{equation} \label{equ:ss}
R^{\s}=R/Rad(R).
\end{equation}
where $Rad(R)$ is the Jacobson radical of $R$. It is well-known that $Rad(R^{\s})=0.$ Additionally since $R$ is Artinian, $R^{\s}$ is Artinian as well. Consequently, $R^{\s}$ is a semisimple algebra. We have the following observation.

\begin{prop}Let $R$ be a finite dimensional $\F_q$ algebra and $Rad(R)$ the Jacobson radical of $R$. Let $\frak{m} \in \mathfrak{P}(R)$. Then 
\begin{enumerate}
    \item $Rad(R) \subset \frak{m} $. 
    \item The map $\frak{m} \mapsto \overline{\frak{m}}:=\frak{m}/Rad(R)$ from $\mathfrak{P}(R) \to \mathfrak{P}(R^{\s})$ is a bijection. Furthermore $N(\frak{m})=N(\frak{\bar{m}}).$
\end{enumerate}

\end{prop}
\begin{proof}
By definition $R/\frak{m} \cong M_r(k)$ for some $r \geq 1$ and a field $k$. The first statement hence follows from \cite[Section 4.3, Lemma b]{pierce1982associative}. The second statement then follows naturally from the first statement.
\end{proof}

A direct consequence of this proposition is the following.
\begin{prop} \label{prop:zeta_semsimplification}
Suppose $R$ is a finite-dimensional $\F_q$-algebra. Then 
\[ \zeta_{R}(s)=\zeta_{R^{\s}}(s) .\] 
\end{prop}

We investigate some further properties of the zeta function of a finite-dimensional $\F_q$-algebra.
\begin{prop}  \label{prop:zeta}
Let $R$ and $T$ be two finite-dimensional $\F_q$-algebras. Then
\begin{enumerate}
    \item $\zeta_{R \times T}(s)= \zeta_{R}(s) \zeta_{T}(s)$. 
    \item If $R$ and $T$ are Morita equivalent, then $\zeta_{R}(s) =\zeta_{T}(s).$
\end{enumerate}

\end{prop}
\begin{proof}
    Part $(1)$ follows directly from the definition of the zeta function. Part $(2)$ is  \cite[Proposition 2.2]{fukaya1998hasse}.
\end{proof}

We discuss some concrete examples of $R$ and their zeta functions. 

\begin{ex}

Let us consider $R=M_n(\F_q).$ Since $M_n(\F_q)$ is Morita equivalent to $\F_q$, Proposition \ref{prop:zeta} shows that 
\[ \zeta_{M_n(\F_q)}(s) = \zeta_{\F_q}(s) = (1-q^{-s})^{-1} .\] 
\end{ex}
\begin{ex}
Let $G$ be a finite group such that $|G|$ is invertible in $\F_q$. Let $R=\F_q[G]$. Suppose further that $G$ splits over $\F_q$; i.e., 
\[ \F_q[G] = \prod_{i=1}^d M_{n_i}(\F_q) .\] 
Then 
\[ \zeta_{\F_q[G]}(s) = \prod_{i=1}^d \zeta_{M_{n_i}(\F_q)}(s)= (1- q^{-s})^{-d} .\] 
\end{ex}
In general, if $G$ does not split over $\F_q$ then the calculation of $\zeta_{\F_q[G]}(s)$ is less explicit. However, when $G$ is abelian, we can explicitly describe the zeta function of $\F_q[G]$. Before we state the key theorem, we recall the following definition. 

\begin{defn}
    Let $d$ be a positive integer and $a$ an integer such that $\gcd(a,d)=1$. The order of $a$ with respect to $d$, denoted by $\ord_d(a)$ is the smallest positive integer $t$ such that $a^t \equiv 1 \pmod{d}.$
\end{defn}

We are now ready to state the key theorem that allows us to compute the zeta function of $\F_q[G]$ where $G$ is an abelian group.

\begin{thm} \cite[Theorem 3.5.4]{milies2002introduction} \label{thm:abelian_factorization}
Let $G$ be a finite abelian group of order $n$ which is prime to $q$. Then 

\[ \F_q[G]  \cong \bigoplus_{d|n} a_d \F_q[\zeta_d] ,\]
where $\zeta_d$ is a primitive root of unity of order $d$ and $a_d=\frac{n_d}{[\F_q(\zeta_d):\F_q]}$. Here $n_d$ is the number of elements of order $d$ in $G.$ Note also that 
\[ [\F_q(\zeta_d):\F_q] = {\rm ord}_{d}(q). \]
\end{thm}
\begin{cor} \label{cor:zeta_abelian_groups}
Let $G$ be a finite abelian group of order $n$ prime to $q$. Then 
\[ \zeta_{\F_q[G]}(s) = \prod_{d|n} (1-q^{-{\rm ord}_{d}(q)s})^{-a_d}, \]
where $a_d$ and ${\rm ord}_{d}(q)$ are as above. 
\end{cor}

We also remark that in some special cases, the zeta function of $\zeta_{\F_q[G]}(s)$ in the modular case (namely $|G|=0$ in $\F_q$) can be deduced from the semisimple case (namely when $|G|$ is invertible in $\F_q$). This is a consequence of Proposition \ref{prop:zeta_semsimplification} and the following theorem. 
\begin{thm} (\cite[Theorem 16.6]{passman1971infinite})
Let $G$ be a finite group. Suppose that $H$ is a normal $p$-Sylow subgroup of $G.$ Then, the Jacobson radical of $\F_q[G]$ is the kernel of the augmentation map 
\[ \epsilon: \F_q[G] \to \F_q[G/H] .\] 
Consequently, $\F_{q}[G]^{\s} \cong \F_q[G/H]$ and 

\[ \zeta_{\F_q[G]}(s) = \zeta_{\F_q[G/H]}(s) .\] 
\end{thm}

We discuss another example of a class of rings where we can compute their zeta functions quite explicitly. Specifically, we can check that the set of all semimagic squares of size $n \times n$ as defined in \cref{defn:semimagic} is a subalgebra of $M_n(k)$. For simplicity, we will denote this ring by $SM_{n}(k)$. By \cite{murase1957semimagic}, we can describe the semisimplification of $SM_n(k)$ explicitly. 

\begin{thm} \cite[Theorem 2, Theorem 3]{murase1957semimagic}
Let $k$ be a field of characteristic $p \geq 0.$ Then 
\begin{enumerate}
    \item If $p \nmid n$ then 
    \[ SM_{n}(k) \cong k \times M_{n-1}(k) .\]
    \item If $p|n$ then the algebra $SM_n(k)$ is not semisimple. Its simplification is given by 
    \[ SM_n(k)^{\s} \cong k \times M_{n-2}(k) .\] 
\end{enumerate}
 
\end{thm}
\begin{cor}
Let $SM_n(\F_q)$ be the ring of all semimagic squares of size $n \times n$ over $\F_q$ with $n \geq 1$. Then 
\begin{enumerate}
\item If $n=1$ then $\zeta_{SM_1(\F_q)}(s)=(1-q^{-s})^{-1}.$
\item If $n=2$ then 
\[ \zeta_{SM_2(\F_q)}(s) = \begin{cases}
  (1-q^{-s})^{-2} &  {\rm char}(\F_q) \neq 2 \\
  (1-q^{-s})^{-1} & {\rm char}(\F_q)=2.
\end{cases} \] 
\item If $n \geq 3$ then \[ \zeta_{SM_n(\F_q)}(s) = (1-q^{-s})^{-2} .\] 

\end{enumerate} 
\end{cor}

We now compute explicitly the zeta function of $\sJ_{G_1, G_2, \ldots, G_d}(\F_q)$ in terms of the zeta functions for $\F_q[G_i]$ for $1 \leq i \leq d.$ We first consider the semisimple case, where all $|G_i|$ are invertible in $\F_q$ (by \cite[Corollary 5.3]{CM2}).  In this case, by Corollary \ref{cor:decomposition}, we have 
\[ \sJ_{G_1, G_2, \ldots, G_d}(\F_q) \cong M_d(\F_q) \times \prod_{i=1}^d \Delta_{\F_q}(G) .\] 
Consequently 
\[ \zeta_{\sJ_{G_1, G_2, \ldots, G_d}(\F_q)}(s) = \zeta_{M_d(\F_q)}(s)  \prod_{i=1}^d \zeta_{\Delta_{\F_q}(G_i)}(s) =(1-q^{-s})^{-1}  \prod_{i=1}^d \zeta_{\Delta_{\F_q}(G_i)}(s) .\] 
Furthermore, we also have 
\[ \F_q[G_i] \cong \F_q \times \Delta_{\F_q}(G_i) ,\]
and therefore 
\[ \zeta_{\F_q[G_i]}(s) = \zeta_{\F_q}(s) \zeta_{\Delta_{\F_q}(G_i)}(s) = (1-q^{-s})^{-1} \zeta_{\Delta_{\F_q}
(G_i)}(s) .\] 
In summary, we have the following 
\begin{prop} \label{prop:zeta_semisimple}
Suppose that $|G_i|$ is invertible in $\F_q$ for $1 \leq i \leq d.$ Then the zeta function of $\sJ_{G_1, G_2, \ldots, G_d}(\F_q)$ is given by 
\[ (1-q^{-s})^{d-1}  \prod_{i=1}^d \zeta_{\F_q[G_i]}(s) .\] 
\end{prop}
We next consider the general case. We can assume that, up to an ordering, there exists a unique positive integer $r$ such that 
\begin{itemize}
    \item $p \nmid |G_i|, 1 \leq i \leq r$. 
    \item $p| |G_i|, r < i \leq d$.
\end{itemize}
We recall the following construction in \cite[Section 5]{CM2}. Let $A$ be a general element of $\sJ_{G_1, G_2, \ldots, G_d}(\F_q)$ 

\begin{equation*}
A=\left[\begin{array}{c|c|c|c}
C_1 & a_{12}\ones & \cdots & a_{1d}\ones \\
\hline
a_{21}\ones & C_2 & \cdots & a_{2d}\ones \\
\hline
\vdots & \vdots & \ddots & \vdots \\
\hline
a_{d1}\ones & a_{d2}\ones & \cdots & C_d
\end{array}\right].
\end{equation*}
We can further partition $A$ into the following blocks 
\[ A =\begin{bmatrix} A_1 & B_1 \\ B_2 & A_2 \end{bmatrix},\] 
where $A_1$ is the union of the upper $r$ blocks, $A_2$ is the union of the lower $d-r$ blocks, $B_1$ (respectively $B_2$) is the union of the upper right (respectively lower left) blocks. Concretely, we have 
\begin{equation*}
A_1=\left[\begin{array}{c|c|c|c}
C_1 & a_{12}\ones & \cdots & a_{1r}\ones \\
\hline
a_{21}\ones & C_2 & \cdots & a_{2r}\ones \\
\hline
\vdots & \vdots & \ddots & \vdots \\
\hline
a_{r1}\ones & a_{r2}\ones & \cdots & C_r
\end{array}\right],
\end{equation*}
\begin{equation*}
A_2=\left[\begin{array}{c|c|c|c}
C_{r+1} & a_{r+1,r+2}\ones & \cdots & a_{r+1,d}\ones \\
\hline
a_{r+2,r+1}\ones & C_{r+2} & \cdots & a_{r+2,d}\ones \\
\hline
\vdots & \vdots & \ddots & \vdots \\
\hline
a_{d,r+1}\ones & a_{d,r+2}\ones & \cdots & C_d
\end{array}\right].
\end{equation*}
Similarly for $B_1, B_2$. Note that we can consider $A_1$ (respectively $A_2$) as an element of $\sJ_{G_1, \ldots, G_r}(\F_q)$ (respectively $\sJ_{G_{r+1}, \ldots, G_d}(\F_q)$.) 

\begin{thm}(\cite{CM2}) \label{thm:jacobson_radical}
Let $I_i$ be the Jacobson radical of $\F_q[G_i]$. Let $\psi$ be the map 
\[ \psi\colon \sJ_{G_1, G_2, \ldots, G_d}(\F_q) \to \sJ_{G_1, \ldots, G_r}(\F_q) \times \prod_{r+1 \leq i \leq d} \F_q[G_i]/I_i ,\]
sending 
\[ A \mapsto (A_1, \overline{C_{r+1}}, \ldots, \overline{C_{d}}) .\] 
Then $\psi$ is a surjective ring homomorphism. Furthermore, the kernel of $\psi$ is the Jacobson radical of $\sJ_{G_1, G_2, \ldots, G_d}(\F_q).$ As a consequence, 
\[ \sJ_{G_1, G_2, \ldots, G_d}(\F_q)^{\s} \cong \sJ_{G_1, \ldots, G_r}(\F_q) \times \prod_{r+1 \leq i \leq d} k[G_i]^{\s} .\] 
Furthermore, by Corollary \ref{cor:decomposition}, we can further decompose 
\[ \sJ_{G_1, G_2, \ldots, G_d}(\F_q)^{\s} \cong M_r(\F_q) \times \prod_{i=1}^r \Delta_{\F_q}(G_i) \times \prod_{r+1 \leq i \leq d} k[G_i]^{\s}. \]
\end{thm}

We have the following corollary by  \cref{thm:jacobson_radical} and  \cref{prop:zeta_semisimple}. 
\begin{thm} \label{cor:general-zeta}
Let $G_1, G_2, \ldots, G_d$ be as above. Then the zeta function of $\sJ_{G_1, G_2, \ldots, G_d}(\F_q)$ is given by 
\[ \zeta_{\sJ_{G_1, G_2, \ldots, G_d}(\F_q)}(s)= (1-q^{-s})^{r-1} \prod_{i=1}^{d} \zeta_{\sJ_{G_i}(\F_q)}(s) .\]  

\end{thm}

\section{$q$-rooted primes and the arithmetic of the join ring $\sJ_{G_1, G_2, \ldots, G_d}(\F_q)$} 
\label{sec:rooted_prime}
In this section, we study the order of the unit group of the join algebra $\sJ_{G_1, G_2, \ldots, G_d}(\F_q)$ where $G_i$ is  cyclic of order $p_i$ with $p_i$ a prime number different from $q$. This is a natural continuation of the work \cite{chebolu2015characterizations} where the author considers the case $d=1$ and $q=2.$ We first recall the following definition.

\begin{defn}
Let $p,q$ be two distinct prime numbers. We say that $p$ is a $q$-rooted prime if $q$ is a primitive root modulo $p$; i.e., $q$ is a generator of the multiplicative group $\F_p^{\times}$ (which is a cyclic group of order $(p-1)$). Equivalently, $p$ is a $q$-rooted prime if and only if $\ord_p(q)=p-1.$
\end{defn}

A conjecture of Emil Artin says that for any non-zero integer $a$ other than $1, -1$ or a perfect square, there exist infinitely many primes $p$ for which $a$ is a primitive root mod $p$. In particular, this would imply that for a given prime $q$, there exist infinitely many $p$ such that $p$ is a $q$-rooted prime. This conjecture remains open, though some partial results are known. For example,
R. Murty and R. Gupta proved unconditionally that there exists an integer that is a primitive root for infinitely many primes.   D. R. Heath-Brown proved that at least one of 2, 3, or 5 is a primitive root modulo infinitely many primes. Furthermore, it is known that Artin's conjecture holds if we assume the Generalized Riemann Hypothesis.  See \cite{gupta1984remark, hooley1967artin} for further discussion on this topic.

In \cite{chebolu2015characterizations}, the authors provide an elegant characterization of $q$-rooted primes using circulant matrices when $q=2$. We remark, however, that their proof remains valid for any prime number $q$. For the sake of completeness, we provide the statement and complete proof here. In the subsequent discussion, the term ``circulant matrices" specifically refers to $G$-circulant matrices where $G$ is a cyclic group.

\begin{thm}  \label{thm:chebolu}
Let $p$ be a prime number. Then, the following statements are equivalent. 
\begin{enumerate}
    \item $p$ is a $q$-rooted prime.
    \item The order of the pole $s=0$ of the zeta function $\zeta_{\F_q[\Z/p]}(s)$ is $2.$
    \item The order of the unit group of the group algebra $\F_q[\Z/p]$ is $(q^{p-1}-1)(q-1).$
    \item The number of invertible circulant matrices of size $p \times p$ over $\F_q$ is $(q-1)(q^{p-1}-1).$
\end{enumerate}
\end{thm}

\begin{proof}
The fact that $(1)$ and $(2)$ are equivalent follows directly from Corollary \ref{cor:zeta_abelian_groups}. For other parts, 
we observe that 
    \[ \F_q[\Z/p] \cong \F_q[x]/(x^p-1) = \F_q \times \F_q[x]/\Phi_p(x).\]
    Here $\Phi_p(x)= \dfrac{x^p-1}{x-1}$ is the $p$-cyclotomic polynomial. By the proof of \cite[Lemma 3.1]{chebolu2015characterizations}, $\Phi_p(x)$ factors as a product of $m=\dfrac{p-1}{\ord_p(q)}$ distinct irreducible polynomials in $\F_q[x]$ of degree $n=\ord_p(q)$. Consequently, as a ring, we have 
    \[ \F_q[\Z/p] \cong \F_q \times \F_{q^n}^{m}. \]
    We see that the order of the unit group of $\F_q[\Z/p]$ is given by $(q-1)(q^n-1)^m.$ We also observe that 
    \[ (q^n-1)^m \leq q^{mn}-1, \]
    and the equality happens iff $m=1.$ This shows the equivalence of (1) and (3). The equivalence of (3) and (4) follows from the observation that units in the group ring correspond to invertible circulant matrices; see Proposition \ref{prop:hurley}.
\end{proof}

The following theorem is a direct generalization of Theorem \ref{thm:chebolu}.
\begin{thm}
Let $q, p_1, p_2, \ldots, p_d$ be prime numbers such that $p_i \neq q.$ Then the following are equivalent
\begin{enumerate}
    \item $p_i$ is a $q$-rooted prime for all $1 \leq i \leq d.$
    \item The order of the pole $s=0$ of the zeta function $\zeta_{\sJ_{\Z/p_1, \Z/p_2, \ldots, \Z/p_d}(\F_q)}$ is $d+1.$
    \item The order of the unit group of the join algebra $\sJ_{\Z/p_1, \Z/p_2, \ldots, \Z/p_d}(\F_q)$ is 
    \[  \prod_{i=1}^d (q^{p_i-1}-1) \times \prod_{i=0}^{d-1} (q^{d}-q^i) .\] 
\end{enumerate}
\end{thm}
\begin{proof}
The equivalence between $(1)$ and $(2)$ follows from \cref{thm:chebolu} and \cref{prop:zeta_semisimple}. Let us consider the other parts. 

By corollary \ref{cor:decomposition} we know that the join algebra $\sJ_{\Z/p_1, \Z/p_2, \ldots, \Z/p_d}(\F_q)$ is decomposed as 
\[ \sJ_{\Z/p_1, \Z/p_2, \ldots, \Z/p_d}(\F_q) \cong M_d(\F_q) \times \prod_{i=1}^d \Delta_{\F_q}(\Z/p_i).\]
Consequently, the order of the unit group of the join algebra $\sJ_{\Z/p_1, \Z/p_2, \ldots, \Z/p_d}(\F_q)$ is given by 
\[ |GL_d(\F_q)| \times \prod_{i=1}^d |\Delta_{\F_q}(\Z/p_i)^{\times}|.\]
By the proof of Theorem \ref{thm:chebolu}, we know that 
\[ |\Delta_{\F_q}(\Z/p_i)^{\times}| \leq q^{p_i-1} -1, \]
with equality when $p_i$ is a $q$-rooted prime.  Combining this with the fact that 
\[ |GL_d(\F_q)| = \prod_{i=0}^{d-1} (q^d-q^i),\]
we get the equivalence of (1) and (3). This completes the proof of the theorem.  
\end{proof}

To motivate another characterization of $q$-rooted primes, consider the following question. \emph{What are all units of order $p$ in the ring $\mathbb{F}_q[{\mathbb{Z}/p}]$, where $p$ and $q$ are distinct primes?} The obvious units of order $p$ that come to mind are of the form $\alpha g$, where $\alpha \in \mathbb{F}_q$ and $g \in \mathbb{Z}/p$ such that $(\alpha g)^p =1$ and $\alpha g \ne 1$. In the literature on group algebras, such units are called trivial units of order $p$. It turns out that these trivial units are the only possible units of order $p$ precisely when $q$ is $p$-rooted.

\begin{prop}
Let $p$ and $q$ be distinct primes. Every unit of order $p$  in $\mathbb{F}_q[\mathbb{Z}/p]$ is trivial if and only if $q$ is $p$-rooted.
\end{prop}

\begin{proof}
Recall the isomorphism:    \[ \F_q[\Z/p] \cong \F_q \times \F_{q^n}^{m}, \]
where $m :=\dfrac{p-1}{\ord_p(q)}$  and  $n :=\ord_p(q)$. Taking units on both sides, we get 
\[ (\F_q[\Z/p])^{\times} \cong \Z/{(q-1)} \times  (\Z/{(q^n-1)})^m.\] 
Every unit of order $p$ in $\mathbb{F}_q[\mathbb{Z}/p]$ is trivial if and only if the number of trivial units of order $p$ in $(\F_q[\Z/p])^{\times}$ is equal to the number of elements of order  $p$ in $\Z/{(q-1)} \times  (\Z/{(q^n-1)})^m$. 

Since the unit group in question is a finite abelian group, the number of elements of order $p$ in $\Z/{(q-1)} \times  (\Z/{(q^n-1)})^m$ is one less than the cardinality of the maximal elementary abelian $p$-subgroup in it. 

Also note that, by definition of $n$, $p$ divides $q^n-1$.  We consider two cases. 
If $p$ does not divide $q-1$, then equating the two  numbers mentioned above, we get 
\[p-1 = p^m-1. \]

Similarly, if $p$ divides $q-1$, we get 
\[p^2-1 = p^{m+1}-1.\]
In both cases, the equations are valid if and only if $m=1$, or equivalently,  $q$ is $p$-rooted.
\end{proof}

We now generalize this to the join of group rings. A unit $u$ in $\sJ_{G_1, G_2, \ldots, G_d}(k)$ is said to be a diagonal unit if $\epsilon(u)$ is a diagonal matrix; that means in our current situation that all off-diagonal blocks of $u$ must be zero. A diagonal unit is trivial if the $i$th  diagonal block is of the form $\circulant([\alpha_i g_i])$, where $\alpha_i \in \F_q$ and $g_i \in G_i$ for $1 \le i \le d$. Note that when $d=1$, this definition gives trivial units for group algebras. We are now ready to state the generalization.

\begin{thm}
Let $q, p_1, p_2, \ldots, p_d$ be prime numbers such that $p_i \neq q.$ Every diagonal unit in $\sJ_{\Z/p_1, \Z/p_2, \ldots, \Z/p_d}(\F_q)$ of order $p$ is trivial if and only if
 $p_i$ is a $q$-rooted prime for all $1 \leq i \leq d$.
\end{thm}

\begin{proof}
Observe that the subgroup of diagonal units in $\sJ_{\Z/p_1, \Z/p_2, \ldots, \Z/p_d}(\F_q)$ is isomorphic to $\prod_{i=1}^d (\F_q[\Z/p_i])^{\times}$. From this, it follows that  every diagonal unit in $\sJ_{\Z/p_1, \Z/p_2, \ldots, \Z/p_d}(\F_q)$ of order $p$ is trivial if and only if every unit of order $p$ is trivial in $\F_q[Z/p_i]$ for each $i$. Invoking the above proposition, we see that the latter holds if and only if $p_i$ is a $q$-rooted prime for all $1 \leq i \leq d$.
\end{proof}

\section{$\sJ_{G_1, G_2, \ldots, G_d}(\F_q)$ and $\Delta_{p^r}$-rings}
\label{sec:delta_ring}
This section considers a special ring-theoretic property of the join ring $\sJ_{G_1, G_2, \ldots, G_d}(\F_q)$. Specifically, we are interested in the $\Delta_n$-property of the join ring. To do so, we first recall the definition of a $\Delta_n$-ring.

\begin{defn} \label{defn:delta_ring}
Let $n$ be a positive integer. A ring $R$ is said to be a $\Delta_n$-ring if for each unit $u \in R^{\times}$, $u^n=1.$ 
\end{defn} 

The $\Delta_n$ property of a ring is well-studied in the literature. It was first introduced in \cite{chebolu2012special}. The author proves that the ring $\Z/n$ of integers modulo $n$ is a $\Delta_2$-ring if and only if $n$ is a divisor of 24.  In \cite{chebolu2013special}, the authors show that the ring $\Z/n[x_1, x_2, \ldots, x_m]$ is a $\Delta_2$-ring  if and only if $n$ is a divisor of 12. Additionally, in \cite{chebolu2015characterizations}, the authors classify all group algebras $k[G]$ which are a $\Delta_p$-ring where $G$ is an abelian group and $p$ is a prime number (see \cite[Theorem  1.4]{chebolu2015characterizations} and \cite[Theorem 1.5]{chebolu2015characterizations}.)

We remark that if $R$ is a $\Delta_n$-ring, it is also a $\Delta_m$ ring if $n|m.$ If $n$ is the smallest positive integer such that $R$ satisfies this property, then call $R$ a strict $\Delta_n$-ring. We will frequently use the fact that whenever a Cartesian product of rings is a $\Delta_n$-ring, so are all the individual factors of the product.
We refer the reader to \cite{chebolu2012special, chebolu2015characterizations, chebolu2013special} for further discussions of this concept.

We next discuss the relationship between the $\Delta_n$-property of $R$ and its semisimplification $R^{\s}.$ For this, we need the following proposition. 

\begin{prop} \label{prop:unit_semi}
    The canonical map $\Phi: R^{\times} \to (R^{\s})^{\times}$ is surjective. 
\end{prop}
\begin{proof} 
    Let $a\in R^{\s}$ be a unit. Then there exists $b\in R^{\s}$ such that $ab=1$. Let $a'\in R$ (respectively, $b'\in R$) be a preimage of $a$ (respectively $b$). One has $a'b'=1+c$ for some $c\in \Rad(R)$. Since $c\in \Rad(R)$, $a'b'=1+c$ is right-invertible. This implies that $a'$ is right-invertible. Similarly, $b'a'=1+d$ for some $d\in \Rad(R)$. From this, we deduce that $a'$ is also left-invertible. Thus $a'$ is a unit and $\Phi(a')=a$. This shows that $\Phi$ is surjective.
\end{proof}

The following lemma follows directly from Proposition \ref{prop:unit_semi}
\begin{lem} \label{lem:semisimplication_delta}
    If $R$ is a $\Delta_n$-ring then so is $R^{\s}.$
\end{lem}

We remark that the converse of Lemma \ref{lem:semisimplication_delta} is not generally true. For example, let $G$ be a $2$-group such as $\Z/8$. Then $R=\F_2[G]$ is a local ring and the Jacobson radical of $R$ is exactly the augmentation ideal $\Delta_{\F_2}(G)$ (see \cite[Corollary 1.4]{carlson2012modules}). Consequently $R^{\s} \cong \F_2$ which is a $\Delta_2$-ring. However, by \cite[Theorem 1.4]{chebolu2015characterizations}, we know that $\F_2[G]$ is not a $\Delta_2$-ring unless $G=(\Z/2)^r.$

When $R$ is a field, we make the following observation that follows immediately from the fact that a polynomial of degree $n$ over a field has at most $n$ roots.

\begin{lem} \label{lem:delta_finite}
    Let $k$ be a field. If $k$ is a $\Delta_n$-ring, then $k$ is a finite field. 
\end{lem}
\begin{rem}
By Lemma \ref{lem:delta_finite}, we can safely assume that all coefficient fields in the discussion below are finite. 
\end{rem}
This section will focus on  $\Delta_{p^r}$-rings where $p$ is a prime number, and $r$ is a positive integer. We remark that the case $r=1$ was studied in \cite{chebolu2015characterizations}, and our work here is a natural continuation of this line of research. We recall that for a group $G$, the exponent of $G$ denoted by $\exp(G)$ is the smallest integer $n$ such that $g^n=1$ for all $g \in G.$ The following simple observation follows directly from the definition of a $\Delta_{p^r}$-ring.

\begin{prop} \label{prop:delta_criterior}
If $R$ is a $\Delta_{p^r}$-ring then $R^{\times}$ is a $p$-group with exponent at most $p^r.$
\end{prop}




Here is an observation that we will use throughout this section.

\begin{lem} \label{lem:matrix_delta}
    Let $q, p^r$ be two prime powers. The matrix algebra $M_n(\F_q)$ is a $\Delta_{p^r}$-ring if and only if $n=1$ and $\F_q$ is a $\Delta_{p^r}$-ring.
\end{lem}
\begin{proof}
Let us assume that $M_n(\F_q)$ is a $\Delta_{p^r}$-ring. By Proposition \ref{prop:delta_criterior}, we know that $GL_n(\F_q)$ must be a $p$-group. Additionally, we know that the order of $GL_n(\F_q)$ is 
\[ \prod_{i=0}^{n-1} (q^n - q^i)=\prod_{i=0}^{n-1}q^{i} (q^{n-i} - 1).\]
Suppose that $n \geq 2.$ We see that $|GL_n(\F_q)| = \prod_{i=0}^{n-1} (q^n - q^i)$ has at least two distinct prime factors. This shows that $GL_n(\F_q)$ is not a $p$-group which contradicts the fact that $M_n(\F_q)$ is a $\Delta_{p^r}$-ring. 
\end{proof}

The main goal of this section is to classify all join algebras $\sJ_{G_1, G_2, \ldots, G_d}(\F_q)$ which are $\Delta_{p^r}$-rings.  To begin this study, we start with the simplest case, namely $\sJ_{G_1, G_2, \ldots, G_d}(\F_q)$ is $\F_q$ (this corresponds to the case $d=1$ and $G_1 = \{e \}$ the trivial group). To answer this question, we first recall the famous Catalan conjecture, now a theorem of Mihailescu (see \cite{mihailescu2004primary, ribenboim1994catalan}.)  

\begin{thm}(See \cite{mihailescu2004primary}) \label{thm:catalant}
The only solution in the natural numbers of the Diophantine equation
 \[ x^{a}-y^{b}=1 .\] 
where $a,b>1$ and $x,y >0$ is $x = 3, a = 2, y = 2, b = 3$.
\end{thm}

Here is a direct corollary of this theorem, which is a generalization of \cite[Lemma 2.1]{chebolu2015characterizations} (see \cite[Theorem 2.4]{chebolu2016fields} for a different but equivalent statement.)

\begin{cor} \label{cor:classification_fields}
    Let $q$ be a prime power. Then the finite field $\F_q$ is a $\Delta_{p^r}$-ring if and only if one of the following conditions hold: 
    \begin{enumerate}
        \item $p=2, q=2^{2^n}+1$ is a Fermat prime, and $r \geq 2^n.$ In this case, $\F_q$ is a strict $\Delta_{2^{2^n}}$-ring.
        \item $p=2^a-1$ is a Mersenne prime and $q=p+1 = 2^a.$  In this case $\F_q$ is a strict $\Delta_p$-ring. 
        \item $p=2$, $q=9$, and $r \geq 3.$ In this case, $\F_q$ is a strict $\Delta_8$-ring.
        \item $q=2$, $p$ and $r$ are arbitrary. 
        
    \end{enumerate}
\end{cor}

\begin{proof}
    The unit group $\F_q^{\times}$ of $\F_q$ is a cyclic group of order $q-1$. Consequently, $\F_q$ is a $\Delta_{p^r}$-ring if and only $q-1|p^r.$ Since $p$ is a prime number, there exists $0 \leq b \leq r$ such that $q-1=p^b.$ Let us write $q=x^a$ where $x$ is a prime and $a$ is a positive integer. We then have the following Diophantine equation 
    \[ x^a - p^b =1.\]

If $a,b>1$, then by Mihailescu's theorem \ref{thm:catalant}, we know that $x=3, a =2, p = 2, b =3,$ thereby satisfying condition (3). If $b=0$ then $x=2, a=1$. Consequently, $q=2$, satisfying condition (4). 
 Next, we consider the case $a=1$. Then $x=p^b+1.$ Since $x>2$, it must be odd. As a result, $p$ is even, hence $p=2$. Therefore, $x=2^b+1.$ From here, we can deduce that $b=2^n$ and $x=2^{2^n}+1$ is a Fermat prime. Thus, in this case, condition (1) is satisfied. Finally, let us consider the case $b=1.$ Then we have $p=x^a-1.$ If $p=2$ then $x=3, a=1$, and we again satisfy condition (1). So, we can safely assume that $p$ is odd. As a result, $x=2$ and $p=2^a-1$ is a Mersenne prime, which fulfills condition (2). 
\end{proof}

Next, we will answer the following question: For which groups $G$ is the group algebra $\F_q[G]$ a $\Delta_{p^r}$-ring? From the canonical embedding $G \hookrightarrow \F_q[G]^{\times}$, we conclude that if $\F_q[G]$ is an $\Delta_{p^r}$-ring then $G$ must be a $p$-group. It turns out that in most cases, $G$ must also be abelian. More precisely, we have the following proposition. 
\begin{prop} \label{prop:abelian}
    Assume that $(p,q) \neq (2, 2)$ and that $\F_q[G]$ is a $\Delta_{p^r}$-ring. Then $G$ is an abelian $p$-group. 
\end{prop}

\begin{proof}
    Since $\F_q \subset \F_q[G]$, we conclude that $\F_q$ is also a $\Delta_{p^r}$-ring. Since $(p,q) \neq (2,2)$, Corollary \ref{cor:classification_fields} implies that $\gcd(p,q)=1.$ Since $G$ is a $p$-group, $|G|$ is invertible in $\F_q.$ By Maschke's theorem, $\F_q[G]$ is semisimple and by the Artin-Wedderburn theorem we must have 
    \[ \F_q[G] \cong \prod_{i=1}^r M_{n_i}(D_i),\]
    where $D_i$ is a division algebra over $\F_q. $ Since $\F_q$ is a finite field, $D_i$ is a finite field as well; see \cite[Chapter 13, Exercise  13, Page 536]{DummitFoote}. By Lemma \ref{lem:matrix_delta}, we conclude that $n_i=1$ and $D_i$ is an $\Delta_{p^r}$-algebra for all $1 \leq i \leq r.$ This implies that $\F_q[G]$ is abelian and hence $G$ is also abelian. 
\end{proof}

We now deal with the case $(p,q)=(2,2)$ separately. Here, instead of working with this particular case, we discuss a more general study of modular group rings, which might be of independent interest. Let $k$ be a finite field of characteristic $p$ and $G$ a finite $p$-group. Let $\Delta_k(G)$ be the augmentation ideal.  It is known that $\Delta_k(G)$ is a nilpotent ideal; in fact $\Delta_k(G)^{|G|}=0$ (see \cite[Corollary 1.3] {carlson2012modules}).  Let $U_1(k[G]):=1 + \Delta_k(G)$ be the set of all normalized units in $k[G].$ We remark that if $u=1+x \in U_1(k[G])$ with $x \in \Delta_k(G)$ then 
\[ u^{|G|}=1+x^{|G|} =1.\]
This shows that $U_1(k[G])$ is a $p$-group. From the isomorphism $k[G]^{\times} \cong k^{\times} \times U_1(k[G])$, we conclude that
$k[G]^{\times}$ is a $p$-group if and only if $k^{\times}$ is a $p$-group. Since $\text{char}(k)=p$, this happens if and only if $k=\F_2.$ In summary, we have
\begin{prop} \label{prop:2_2_case}
    Let $(p,q)=(2,2)$ and $G$ a $2$-group. Then $\F_q[G]$ is a $\Delta_{2^r}$-ring where $2^r=\exp(U_1(\F_q[G]))$. Furthermore, if $G$ is abelian, $\F_q[G]$ is a strict $\Delta_{\exp(G)}$-ring. 
\end{prop}
\begin{proof}
    We already explained the proof of the first part. For the second part, we note that if $G$ is abelian and 
    \[ u = 1 +x = \sum_{g \in G} a_g g \]
    is a normalized unit (so $a_e=1$), then 
    \[ u^{\exp(G)} = \sum_{g \in G} a_g^{\exp(G)} g^{\exp(G)} = \sum_{g \in G} a_g = 1.\]
\end{proof}

\begin{rem}
It is worth mentioning that the problem of determining $\exp(U_1(\F_2[G])$ is well-studied but still open in the literature. For further discussions, see \cite{johnson1978modular, shalev1991lie}. 
\end{rem}

With these preliminary results, we are now ready to classify all group algebras $\F_q[G]$ which are $\Delta_{p^r}$-rings.

\begin{thm} \label{thm:delta_classification}
    Let $q, p^r$ be prime powers. Let $G$ be a finite group. The group algebra $\F_q[G]$ is a $\Delta_{p^r}$-ring if and only if $G$ is a $p$-group and one of the following conditions holds. 
    \begin{enumerate}
        \item $p=2, q=2^{2^n}+1$ is a Fermat prime, $r \geq 2^n$, $G$ is abelian, and the exponent of $G$ is a divisor of $2^{2^n}.$ In this case, $\F_q[G]$ is a strict $\Delta_{2^{2^n}}$-ring. 
        \item $p=2^a-1$ is a Mersenne prime, $q=p+1 = 2^a$ and $G = (\Z/p)^s$ for some $s \geq 0.$ In this case $\F_q[G]$ is a strict $\Delta_{p}$-ring. 
        \item $p=2^a-1$ is a Mersenne prime, $q=2$ and $G=(\Z/p)^s$ for some $s \geq 0.$ 
        \item $p=2, q=3$, $r\geq 3$, $G$ is abelian, and the exponent of $G$ is $4$ or $8.$
        \item $p=2, q=9, r \geq 3$, $G$ is abelian, and the exponent of $G$ is at most $8.$ In this case, $\F_q[G]$ is a strict $\Delta_8$-ring.
        \item $q=2, p=2$ and $2^r \geq \exp(U_1(\F_2[G]))$. 
        
    \end{enumerate}
\end{thm}
\begin{proof} 



We will discuss both directions of the above theorem simultaneously. First, note that since $\F_q$ is a subring of $\F_q[G]$, if $\F_q[G]$ is a $\Delta_{p^r}$-ring then so is $\F_q.$  On the other hand, from Corollary \ref{cor:classification_fields}, $\F_q[G]$ is a $\Delta_{p^r}$-ring if any of conditions (1)-(4) hold.  Thus, we may assume that  $\F_q[G]$ is a $\Delta_{p^r}$-ring (since otherwise, neither side of our equivalence holds).  Also from Corollary \ref{cor:classification_fields}, $p \neq \text{char}(\F_q)$ unless $p=q=2.$ The case $(p,q)=(2,2)$ is treated separately in Proposition \ref{prop:2_2_case}. For now, let us assume that $(p,q) \neq (2,2)$. 
By Proposition \ref{prop:abelian}, we conclude that $G$ is abelian.  Since $G$ is an abelian $p$-group with $\gcd(p,q)=1$, Theorem \ref{thm:abelian_factorization} implies
\[ \F_q[G]  \cong \bigoplus_{d||G|} a_d \F_q[\zeta_d].\]
Here $\zeta_d$ is a primitive root of unity of order $d$ and $a_d=\frac{n_d}{[\F_q(\zeta_d):\F_q]}$ where $n_d$ is the number of elements of order $d$ in $G.$ From this formula, we conclude that $\F_q[G]$ is a $\Delta_{p^r}$-ring if and only if each component $\F_q[\zeta_d]$ is. Since $|G|$ is a $p$-group and $\F_q[\zeta_{d'}] \subset \F_q[\zeta_{d}]$ if $d'|d$, we conclude that $\F_q[G]$ is a $\Delta_{p^r}$-ring if and only if $\F_q[\zeta_D]$ is $\Delta_{p^r}$-ring where $D$ is largest number such that $n_D>0.$ Since $G$ is a $p$-group, $D$ is exactly the exponent of $G.$ We remark that $\F_q[\zeta_D] = \F_{q^m}$ where 
\[ m = [\F_q(\zeta_D):\F_q] = \text{ord}_{D}(q). \]

We now consider a few cases based on the classification described in Corollary \ref{cor:classification_fields}. 

\textbf{Case 1}: $p=2$ and $q^m = 2^{2^n}+1$ is a Fermat prime. This shows that $m=1$ and $q$ is a Fermat prime. Furthermore, $m=1$ means that $\text{\ord}_{D}(q)=1$ or equivalently $D|q-1 = 2^{2^n}.$ This covers the first case of our theorem. 

\textbf{Case 2}: $p=2^a-1$ is a Mersenne prime and $q^m=p+1=2^a.$ This shows that $q=2^b$ with $bm = a.$ By definition of $m$, we have
$q^m \equiv 1 \pmod{D}$.  Since $q^m=2^a$, this is equivalent to $D|2^a-1=p.$ This implies $D=1$ or $D=p.$ From this, we can conclude that $G=(\Z/p)^s$ for some $s \geq 0.$ Furthermore, we remark that since $p=2^a-1$ is a prime number, $a$ is a prime number. We then see that $(b,m)=(a,1)$ or $(b,m)=(1,a).$ The case $(b,m)=(a,1)$ covers the second case of our theorem and the case $(b,m)=(1,a)$ covers the third case of our theorem. 

\textbf{Case 3:} $p=2$, $q^m=9$ and $r\geq 3$. First, consider the case where $(q,m)=(3,2).$ Since $m=2$, we know that $9=q^m \equiv 1 \pmod{D}$ and $q \not \equiv 1 \pmod{D}.$ This shows that $D \in \{4, 8 \}.$ This covers the fourth case of our theorem. Next, consider the case $(q,m)=(9,1).$ Again, we see that $D|8.$ This covers the fifth case of our theorem. 
  \end{proof}

We now focus on a special case. 
\begin{defn} \label{defn:diagonal_property}
A ring $R$ is said to have the \textit{diagonal property} if it is a $\Delta_2$-ring. 
\end{defn}

The classification given in Theorem \ref{thm:delta_classification} provides another proof for the following statements, which were first proved in \cite{chebolu2015characterizations} under the assumption that $G$ is abelian.
\begin{cor}(\cite[Theorem 1.4]{chebolu2015characterizations} and \cite[Theorem 1.5]{chebolu2015characterizations})  Let $G$ be a group and $k$ a field.\label{cor:chebolu_classification}
\begin{enumerate}

\item  The group algebra $k[G]$ has the diagonal property if and only if $k[G]$ is either $\F_2[(\Z/2)^r]$ or $\F_3[(\Z/2)^r].$
\item  Let p be an odd prime. The group algebra $k[G]$ is a $\Delta_p$-ring if and only if $p$ is a Mersenne prime and $k[G]$ is is either $\F_2[(\Z/p)^r]$ or $\F_{p+1}[(\Z/p)^r].$
\end{enumerate}

\end{cor}

Thus, these results give us a simple and elegant characterization of Mersenne primes.  An odd prime $p$ is a Mersenne if and only if $\F_2[E]$ is a $\Delta_p$-ring, where $E$ is any finite elementary abelian $p$-group.

Finally, we answer the following question: which join algebra $\sJ_{G_1, G_2, \ldots, G_d}(\F_q)$ is a $\Delta_{p^r}$-ring. 
\begin{thm}
Suppose that $d \geq 2.$ Then the join algebra $\sJ_{G_1, G_2, \ldots, G_d}(\F_q)$ is a $\Delta_{p^r}$-ring if and only if the following conditions are satisfied
\begin{enumerate}
    \item $p=q=2.$
    \item $G_i$ is a $2$-group for all $1 \leq i \leq d$.
    \item There is at most one index $i$ such that $G_i = \{e \}$ the trivial group. 
    \item $2^r \geq \max_{1 \leq i \leq d} \exp(U_1(\F_2[G_i]))$.
\end{enumerate}
\end{thm}
\begin{proof}
Let us prove the ``only if'' part of the above theorem. So assume $\sJ_{G_1, G_2, \ldots, G_d}(\F_q)$ is a $\Delta_{p^r}$-ring.  First, we claim that $(p,q)=2.$ In fact, suppose that $(p,q) \neq (2,2)$. Let us consider the following embedding $\F_q[G_d]^{\times} \hookrightarrow \sJ_{G_1, G_2, \ldots, G_d}(\F_q)^{\times}$ sending  

\[ C_1 \mapsto \left[\begin{array}{c|c|c|c}
C_1 &  0  & \cdots & 0  \\
\hline
0 & I_{k_2} & \cdots & 0  \\
\hline
\vdots & \vdots & \ddots & \vdots \\
\hline
0 & 0 & \cdots & I_{k_d}
\end{array}\right].
\]
Since $\sJ_{G_1, G_2, \ldots, G_d}(\F_q)$ is a $\Delta_{p^r}$-ring, $k[G_1]$ is a $\Delta_{p^r}$-ring as well. Similarly, $k[G_i]$ is a $\Delta_{p^r}$-ring for all $ 1 \leq i \leq d$. We conclude that $G_i$ is a $p$-group for all $1 \leq i \leq d$. Furthermore, by Theorem \ref{thm:delta_classification}, we know that $\gcd(p,q)=1$ since we assume that $(p,q) \neq (2,2)$. Then by Corollary \ref{cor:decomposition}, $M_d(\F_q)$ is a direct factor of $\sJ_{G_1, G_2, \ldots, G_d}(\F_q)$. This shows that $M_d(\F_q)$ is a $\Delta_{p^r}$-ring as well.  However, Lemma \ref{lem:matrix_delta} implies that $d=1$, which is a contradiction. This shows that $(p,q)=(2,2).$

From now on, we will assume that $(p,q)=2.$ In particular, this implies that $G_i$ is a $2$-group. Suppose there are exactly $t$ elements amongst $G_i$, which are trivial groups. We claim that $t \leq 1$. In fact, by Theorem \ref{thm:jacobson_radical}, $M_t(\F_2)$ is a direct factor of $\sJ_{G_1, G_2, \ldots, G_d}(\F_2)^{\s}$ which is a $\Delta_{p^r}$-ring by Lemma \ref{lem:semisimplication_delta}.  This shows that $M_t(\F_2)$ is a $\Delta_{p^r}$-ring. By Lemma \ref{lem:matrix_delta}, we conclude that $0 \leq t \leq 1.$ Finally, the embedding $\F_2[G_i]^{\times} \hookrightarrow \sJ_{G_1, G_2, \ldots, G_d}(\F_2)^{\times}$ explained above implies that 
\[ 2^r \geq  \exp(U_1(\F_2[G_i]), \forall 1 \leq i \leq d .\] 

In summary, we have proved the ``only if" part of the theorem. We now prove the converse. Let us consider the case that all $G_i$ are nontrivial $2$-groups. Let 
\begin{equation*}
A=\left[\begin{array}{c|c|c|c}
C_1 & a_{12}\ones & \cdots & a_{1d}\ones \\
\hline
a_{21}\ones & C_2 & \cdots & a_{2d}\ones \\
\hline
\vdots & \vdots & \ddots & \vdots \\
\hline
a_{d1}\ones & a_{d2}\ones & \cdots & C_d
\end{array}\right].
\end{equation*}
be an invertible element in $\sJ_{G_1, G_2, \ldots, G_d}(\F_2)$. Then $\epsilon(A)$ is invertible where $\epsilon: \sJ_{G_1, G_2, \ldots, G_d}(\F_q) \to M_d(\F_2)$ is the augmentation map. By definition, we have 
\begin{equation*}
\epsilon(A)=\left[\begin{array}{c|c|c|c}
\epsilon(C_1) &  0  & \cdots & 0  \\
\hline
0 & \epsilon(C_2) & \cdots & 0  \\
\hline
\vdots & \vdots & \ddots & \vdots \\
\hline
0 & 0 & \cdots & \epsilon(C_d)
\end{array}\right] \in GL_d(\F_2).
\end{equation*}
We conclude that $\epsilon(C_i)=1$ for all $1 \leq i \leq d.$ This implies that $C_i$ is invertible for $1 \leq i \leq d$ since $\F_2[G_i]$ is a local ring in which $\Delta_{\F_2}(G_i)$ is the maximal ideal. We then see that 

\begin{equation*}
A^2=\left[\begin{array}{c|c|c|c}
C_1^2 &  0  & \cdots & 0  \\
\hline
0 & C_2^2 & \cdots & 0  \\
\hline
\vdots & \vdots & \ddots & \vdots \\
\hline
0 & 0 & \cdots & C_d^2
\end{array}\right].
\end{equation*}

Consequently 
\begin{equation*}
A^{2^r}=\left[\begin{array}{c|c|c|c}
C_1^{2^r} &  0  & \cdots & 0  \\
\hline
0 & C_2^{2^r} & \cdots & 0  \\
\hline
\vdots & \vdots & \ddots & \vdots \\
\hline
0 & 0 & \cdots & C_d^{2^r}
\end{array}\right].
\end{equation*}

Since $\F_2[G_i]$ is a $\Delta_{2^r}$-ring, we conclude that $C_i^{2^r} = I.$ As a result, $A^{2^r}=I.$ This shows that $\sJ_{G_1,G_2, \ldots, G_d}(\F_2)$ is a $\Delta_{2^r}$-ring.  The case where there is one $G_i = \{e \}$ can be proved using similar calculations. 
\end{proof}

A direct corollary of the above theorem is the following.
\begin{cor} \label{cor:join_algebra_diagonal}
Suppose that $d \geq 2.$ Then the join algebra $\sJ_{G_1, G_2, \ldots, G_d}(\F_q)$ has the diagonal property if and only if $\sJ_{G_1, G_2, \ldots, G_d}(\F_q)$ is $\sJ_{(\Z/2)^{r_1}, (\Z/2)^{r_2}, \ldots, (\Z/2)^{r_d}}(\F_2)$ where $r_i \in \Z_{\geq 0}$ and at most one of the $r_i$ is equal to $0.$

\end{cor}

\bibliographystyle{abbrv}
\bibliography{CM2.bib}

\end{document}